\documentclass[10pt,conference]{IEEEtran}
\IEEEoverridecommandlockouts
\usepackage{geometry}
\usepackage{afterpage}
\usepackage{cite}
\usepackage{graphicx}
\usepackage{subfigure}
\usepackage{booktabs}
\usepackage{tabularx}
\usepackage{multirow}
\usepackage[font=small]{caption}
\usepackage[hyphenbreaks]{breakurl}
\usepackage{textcomp}
\usepackage{xcolor}
\usepackage{diagbox}
\usepackage{float}
\usepackage{amsmath,amssymb,amsfonts}
\usepackage{hyperref}
\usepackage[hyphenbreaks]{breakurl}
\usepackage{url}
\usepackage{datatool}
\usepackage{algorithm}
\usepackage[noend]{algpseudocode}
\usepackage{enumitem}

\makeatletter

\newcommand{\Rmnum}[1]{\expandafter\@slowromancap\romannumeral #1@}
\makeatother

\usepackage[textsize=tiny, textwidth=1.5cm]{todonotes}
\usepackage{xcolor, soul}
\sethlcolor{yellow}

\usepackage{listings}
\AtBeginDocument{%
  \providecommand\BibTeX{{%
    \normalfont B\kern-0.5em{\scshape i\kern-0.25em b}\kern-0.8em\TeX}}}

\begin{document}

\title{Grid-Aware Charging and Operational Optimization for Mixed-Fleet Public Transit}

\author{
    \IEEEauthorblockN{Rishav Sen\IEEEauthorrefmark{1}, Amutheezan Sivagnanam\IEEEauthorrefmark{2}, Aron Laszka\IEEEauthorrefmark{2}, Ayan Mukhopadhyay\IEEEauthorrefmark{1}, Abhishek Dubey\IEEEauthorrefmark{1}}
    \IEEEauthorblockA{\IEEEauthorrefmark{1} Institute for Software Integrated Systems, Vanderbilt University}
    \IEEEauthorblockA{\IEEEauthorrefmark{2} Pennsylvania State University}
}

\maketitle

\begin{abstract}
The rapid growth of urban populations and the increasing need for sustainable transportation solutions have prompted a shift towards electric buses in public transit systems. However, the effective management of mixed fleets consisting of both electric and diesel buses poses significant operational challenges. One major challenge is coping with dynamic electricity pricing, where charging costs vary throughout the day. Transit agencies must optimize charging assignments in response to such dynamism while accounting for secondary considerations such as seating constraints. This paper presents a comprehensive mixed-integer linear programming (MILP) model to address these challenges by jointly optimizing charging schedules and trip assignments for mixed (electric and diesel bus) fleets while considering factors such as dynamic electricity pricing, vehicle capacity, and route constraints. We address the potential computational intractability of the MILP formulation, which can arise even with relatively small fleets, by employing a hierarchical approach tailored to the fleet composition. By using real-world data from the city of Chattanooga, Tennessee, USA, we show that our approach can result in significant savings in the operating costs of the mixed transit fleets. 
\end{abstract}

\begin{IEEEkeywords}
Mixed transit fleet, electrification, dynamic pricing, hierarchical MILP
\end{IEEEkeywords}

\section{Introduction}

As global concerns about climate change and sustainable energy continue to grow, the transition to electric buses in public transportation systems has become a crucial strategy in reducing greenhouse gas emissions and promoting sustainable urban mobility \cite{mckinsey2021integrated, litman2010evaluating}. 
However, the transition has been slow due to the vast difference in cost between electric and diesel-powered buses, e.g., while a diesel transit bus costs about \$500,000, electric buses cost more than \$750,000~\cite{buscost}. As a result, most transit agencies operate mixed fleets, which allows them to slowly supplement their existing fleet with electric vehicles, leading to a fleet that has both diesel and electric vehicles and may also consist of hybrid vehicles.

As the transit agencies increase the ratio of electric buses in their fleet, they face challenges not only in deciding which bus should be used on what trips but also in managing when the EVs charge; this problem is particularly taxing due to the comparatively shorter range of electric buses, which is affected by weather and traffic. Moreover, optimizing such fleets must be done while ensuring that routes are served reliably without delays and accommodating all commuters, such as having seating requirements. 
Another unique challenge faced by electric fleets is dynamic electricity pricing~\cite{amin2020review}. The cost of purchasing electricity varies over time for commercial customers such as transit agencies, with significant variability within each day. Dynamic electricity pricing usually consists of a Time-of-use (TOU) based price \cite{widergren2012real}. The usage charge typically varies over smaller temporal granularity such as a day, while the demand charge is aggregated over a longer duration such as a month. As a result, \textit{ad-hoc} charging strategies can be costly, thereby making it imperative for agencies to plan ahead. Diesel prices exhibit significantly lesser variance within a day, and transit agencies do not incur any demand-based charges for the fuel. 

The problem of optimal charging assignment for electric vehicles has been studied, e.g., there is extensive prior work on vehicle assignment to trips by considering capacities and availability of charging infrastructure~\cite{sivagnanam2021minimizing, hu2016electric, cooney2013life, buekers2014health}. However, to the best of our knowledge, there is no prior work on the significantly more complex problem of integrating all the nuances of mixed-fleet transit, such as dynamic electricity pricing and seating constraints, while optimizing traditional transit objectives such as assigning mixed fleets to trips on fixed lines. We reiterate that the assignment problem is itself challenging, given the concerns about the driving range and battery capacity of EVs. In this paper, we present a problem formulation that explores these unique challenges in collaboration with a public transit agency of a medium-sized city in the U.S.A. Our partner agency operates both electric and diesel buses, serving about 200,000 residents.

\noindent \textbf{Contributions:} Our formulation is inspired by prior work by Sivagnanam et al.~\cite{sivagnanam2021minimizing}, but we make three major contributions. 

\begin{enumerate}[leftmargin=*]
    \item We extend the formulation from \cite{sivagnanam2021minimizing} to enable us to focus on both the objectives of energy consumption as well as the efficient trip assignment while adhering to the operational constraints of charging and vehicle range.
    Our formulation balances the twin objectives of energy consumption and efficient trip assignment while adhering to the operational constraints of charging and vehicle range.
    
    \item We incorporate seating restrictions that may apply on a given transit block as the transit agency might have requirements of assigning specific capacity buses on certain blocks, due to local rules, or known mobility issues with the people being served by the buses.
    
   \item  Existing exact solution approaches are intractable in practice. Heuristic solutions can scale, but they can sacrifice performance. We demonstrate that unlike an approximated simulated annealing approach as described in \cite{sivagnanam2021minimizing}, we can solve the problem with better results using a hierarchical approach where we iteratively solve related Mixed Integer Linear Programming (MILP) problems, describing the decomposition technique and other computational approaches to make it more tractable (e.g. the warm start to the solution).

\end{enumerate}

 Our approach aims to provide a near-optimal solution while maintaining the feasibility of the solution. This approach not only enhances the practicality of fleet management strategies but also paves the way for more sustainable and economically feasible public transportation systems. Using our method, we show an average improvement of $2.58\%$ in operational cost over the foundational work of Sivagnanam et al. \cite{sivagnanam2021minimizing}, and an average $6.25\%$ improvement over the actual operational cost that the transit agency currently incurs.

\section{Related Research}\label{sec:related}



The shift towards electrification in public transit is driven by the urgent need to reduce greenhouse gas emissions and enhance urban air quality. Studies have extensively compared the environmental impact and lifecycle costs of electric versus diesel buses, highlighting the long-term benefits despite higher initial costs \cite{quarles2020costs, potkany2018comparison} and underlining the sustainability advantages of electric buses, helping meet climate goals. However, the operational challenges posed by electric buses remain a significant concern for transit agencies, showing a need for enhanced management strategies to ensure a smoother transition.

Optimization of fleet operations has been a critical focus area, with several studies on efficiently managing mixed fleets of electric and diesel buses. There have been discussions on optimization approaches for vehicle assignment considering the availability of charging infrastructure and the constraints imposed by battery limits \cite{sivagnanam2021minimizing, mohammadbagher2022multi}. Similarly, Hu et al. \cite{hu2016electric} explore algorithms that optimize the scheduling of electric buses to maximize operational efficiency while minimizing charging times and costs. These models are foundational, however, they often overlook the complexities introduced by dynamic electricity pricing, significantly altering the cost-benefit landscape.

Dynamic electricity pricing presents a unique challenge for electric bus operations. The variability in electricity costs due to Time-of-use (TOU) can greatly influence charging strategies. Hao et al. provide a comprehensive review of how dynamic pricing affects large-scale consumers like transit agencies, emphasizing the need for advanced planning and real-time decision-making frameworks \cite{hao2024dynamic}. It has also been studied in the context of personal electric vehicle charging \cite{soares2017dynamic}. Despite the insights offered by these studies, there is a notable gap in research regarding the integration of dynamic pricing models into the operational strategies of mixed-fleet transit systems, which our research aims to address.

Hierarchical optimization techniques, particularly hierarchical MILP, have been successfully applied in various industrial and systems engineering contexts to decompose complex decision-making problems into more manageable sub-problems. This method allows for detailed planning at different decision levels, facilitating nuanced control over complex systems with multiple interacting components. For instance, in the context of logistics and supply chain management, hierarchical MILP has been used to optimize resource allocation and electrical load management \cite{zhang2021hierarchical} \cite{nebuloni2023hierarchical}, cooperative vehicle networks \cite{1657382}. Although hierarchical approaches have been hinted at in transportation literature \cite{yokoyama2021hierarchical}, their application has predominantly been limited to simpler scenarios that do not account for the dynamic and multi-faceted challenges presented by mixed-fleet transit systems under variable pricing conditions.

Existing literature extensively explores separate elements of fleet optimization, dynamic pricing impacts, and the potential of hierarchical MILP approaches. However, it often overlooks the integration of these elements into a cohesive model tailored for mixed-fleet transit operations. Our research bridges this gap by adapting hierarchical MILP to handle the intricacies of mixed-fleet management under dynamic pricing conditions, incorporating operational constraints like vehicle range, battery limitations, and seating capacities. This novel application extends previous optimization frameworks and introduces a multi-layered decision-making structure that aligns with the operational realities of modern transit agencies.
\section{Problem Model}\label{sec:model}

\noindent \textbf{Vehicles} Our area of consideration is a transit agency that operates a set of buses $\mathcal{V}$, where each bus $v \in \mathcal{V}$ belongs to a vehicle model $\mathcal{M}_{v} \in \mathcal{M}$ = $\mathcal{M}^{\text {diesel}} \cup \mathcal{M}^{\text {electric }}$. Each vehicle model has a state of charge (SoC) between$\left [{\mathcal{M}_{v}}^{min}, {\mathcal{M}_{v}}^{max} \right ]$. Each vehicle also has an operating efficiency $v^{op}$, the average energy used per mile. We consider hybrid buses under the set $\mathcal{M}^{\text{diesel}}$ as they use diesel fuel, and do not require charging. Hybrid buses have higher $v^{op}$ than regular diesel buses, and electric buses have the highest efficiency. We use an operating efficiency threshold $v^{op}_{th}$ to form the hierarchical formulation. Each vehicle has a seating capacity of $v^s$ seats.

\noindent \textbf{Locations} Locations $\mathcal{L}$ include terminals (with charging stations), bus stops, and additional charging stations in the transit network. 

\noindent \textbf{Transit blocks} are a fundamental concept in transit operations management, encompassing a set of sequential trips assigned to a single vehicle for a day of operation. The use of transit blocks helps optimize vehicle and driver scheduling to operate efficiently and reliably, by minimizing vehicle downtime, reducing operational costs, and streamlining the planning process. A bus serving a block $t \in \mathcal{T}$ leaves from the block origin $t^{\text {origin }} \in \mathcal{L}$ at time $t^{\text {start }}$ and arrives at destination $t^{\text {destination }} \in \mathcal{L}$ at time $t^{\text {end}}$. The distance it covers is $t^d$.  Each block may require the buses to have a minimum number of seats, denoted by $t^s$.

\noindent \textbf{Charging rate} It is the amount of power delivered to the electric buses, measured in terms of kilowatt-hour (kWh). Each charger has a specified minimum charging limit of $0$ kWh and a maximum charging limit $q^{max}$ kWh. We assume all chargers here are unidirectional, but this formulation can be easily extended to the use of bidirectional chargers, changing the minimum to be negative ($q^{min}$ kWh, $q^{min} < 0$).

\noindent \textbf{Time slots} We denote uniform-length time slots using $\mathcal{S}$, which are used to divide the operational hours into different discrete intervals and are used for both block and charging assignment. A time slot $s \in \mathcal{S}$ begins at $s^{\text {start }}$ and ends at $s^{\text {end }}$. The very last time slot of the day is denoted by $s_{last}$ and is used to replenish the vehicles to their maximum capacity. We assume this is the downtime after a day's operation, having no more blocks to serve and this time slot can continue till the start of the next day's operations.

\noindent \textbf{Charging} We denote the set of charger poles $\mathcal{C P}$, where $c p^{\text {location }} \in \mathcal{L}$ is the location of charging pole $cp \in \mathcal{CP}$.  Charging times can start and end only at the start and end of each time slot. A charging pole $cp \in \mathcal{CP}$ can provide up to $q_{cp}^{max}$ amount of energy to an electric bus $v$ at time slot $s$, denoted by $c^v_s$. We assume this amount of charge is uniformly fed to the bus over the entire duration of the slot. Let $\mathcal{C}=\mathcal{C} \mathcal{P} \times \mathcal{S}$ be the set of all charging slots. We assume that the buses are fully charged and fuelled at the start of the day, and at the end of the day,  we replenish them back to full capacity.

\noindent \textbf{Non-service trips} Besides serving transit blocks, buses may also need to drive between blocks or charging poles. For example, if a bus has to serve a block that starts from a location different from the previous block's destination, the bus first needs to drive to the origin of the subsequent block. An electric bus may also need to drive to a charging pole after serving a transit block to recharge, then drive from the pole to the origin of the next transit block. We will refer to these trips, which are driven outside of revenue service, as non-service trips. $T\left(l_{1}, l_{2}\right)$ denotes the non-service trip from location $l_{1} \in \mathcal{L}$ to $l_{2} \in \mathcal{L}$; and \textcolor{red} {$D\left(l_{1}, l_{2}\right)$} denotes the time duration of this non-service trip.

\noindent \textbf{Energy usage} The energy used for every transit block or non-service trip is denoted by $E(v,t) \ \forall t \in \mathcal{T}, v \in \mathcal{V}$. 

\noindent \textbf{Distance traveled} The distance traveled by a vehicle for a non-service trip between blocks or chargers $x_1, x_2$ is represented using \textcolor{red}{$\Delta(x_1, x_2)$}.

\noindent \textbf{Time-of-Use electricity price} We opt to use Time-of-Use pricing to formulate this model, with multiple pricing tiers, $\{w_1, w_2, \cdots w_{n} \} \in \mathcal{W}$ where $w_s$ is the electricity price per kWh during the time slot, $s$.

\noindent \textbf{Operating cost} We need to recharge electric buses and refuel the diesel buses to replenish them. We use the common herm \textit{replenish} to refer to both \textit{recharging} and \textit{refueling} of electric and diesel use respectively. The variable $g_s$ represents the cost of replenishing a bus during time slot $s$. $\sum_{s} g_s$ gives the total operating cost.

\subsection{Solution Space}\label{se:soln_space}

\textbf{Assignments}
$\mathcal{A}$ is the set of solutions, where for each block $t \in \mathcal{T}$, exactly one bus $v \in \mathcal{V}$ is assigned to serve block $t$, such that $\langle v, t \rangle \in \mathcal{A}$. Also, each electric bus $v$ must be charged before its battery state of charge drops below the safe level for operation, ${\mathcal{M}_{v}}^{min}$. At most one electric bus $v$ can be assigned to one charging slot $(c p, s) \in \mathcal{C}$ and is denoted as $\langle v,(c p, s)\rangle \in \mathcal{A}$. We consistently use the assumption that an electric bus remains at the charging pole for the entire duration of the charging slot, and is charged at a uniform rate for the entire time slot.

\vspace{2mm}


\noindent \textbf{Feasibility Constraints}\label{sec:feasibility}
If a bus $v$ is assigned to serve an earlier transit block $t_{1}$ and a later block $t_{2}$, then the duration of the non-service trip from $t_{1}^{\text {destination }}$ to $t_{2}^{\text {origin }}$ must be less than or equal to the time between $t_{1}^{\text {end }}$ and $t_{2}^{\text {start }}$. Otherwise, it would not be possible to serve $t_{2}$ on time. We formulate this constraint as:
\small
\begin{multline*}
\forall t_{1}, t_{2} \in \mathcal{T} ; t_{1}^{\text {start }} \leq t_{2}^{\text {start }} ;\left\langle v, t_{1}\right\rangle \in \mathcal{A} ;\left\langle v, t_{2}\right\rangle \in \mathcal{A}: \\
t_{1}^{\text {end }}+D\left(t_{1}^{\text {destination }}, t_{2}^{\text {origin }}\right) \leq t_{2}^{\text {start }}
\end{multline*}
\normalsize
If the constraint is satisfied by every pair of consecutive blocks assigned to a bus, then it is also satisfied by every pair of non-consecutive blocks assigned to the bus.
We need to formulate similar constraints for non-service trips to, from, and between charging slots:
\small
\begin{multline*}
\forall t \in \mathcal{T} ;(c p, s, q) \in \mathcal{C} ; t^{\text {start }} \leq s^{\text {start }} ;\langle v, t\rangle,\langle v,(c p, s, q)\rangle \in \mathcal{A}: \\
t^{\text {end }}+D\left(t^{\text {destination }}, c p^{\text {location }}\right) \leq s^{\text {start }} \\
\end{multline*}
\vspace{-13mm}
\begin{multline*}
\forall t \in \mathcal{T} ;(c p, s, q) \in \mathcal{C} ; t^{\text {start }} \geq s^{\text {start }} ;\langle v, t\rangle,\langle v,(c p, s, q)\rangle \in \mathcal{A}: \\
s^{\text {end }}+D\left(c p^{\text {location }}, t^{\text {origin }}\right) \leq t^{\text {start }}  \\
\end{multline*}
\vspace{-13mm}
\begin{multline*}
\forall\left(c p_{1}, s_{1}, q_1\right),\left(c p_{2}, s_{2}, q_2\right) \in \mathcal{C} ; s^{\text {start }} \leq s^{\text {start }} ;\\
\left\langle v,\left(c p_{1}, s_{1}, q_1\right)\right\rangle,\left\langle v,\left(c p_{2}, s_{2}, {q_2}\right)\right\rangle \in \mathcal{A}: \\
s_{1}^{\text {end }}+D\left(c p_{1}^{\text {location }}, c p_{2}^{\text {location }}\right) \leq s_{2}^{\text {start }} 
\end{multline*}
\normalsize
The above four equations are collectively represented using $F(x_1, x_2)$, which form the feasibility checks.

\subsection{Objective}

This objective equation represents a minimization problem in which the goal is to minimize the total cost of replenishing the buses.
\small
\begin{equation}
\min \sum_{s\in S} g_s 
\end{equation}
\normalsize
\noindent The equation consists of two main components:

Recall $\sum_{s\in S} g_s$, represents the total cost of recharging and refueling the buses.
By minimizing this sum, we can find the optimal bus assignment that minimizes the overall cost of replenishing the buses.

\section{Hierarchical Approach}\label{sec:exact_soln} \label{sec:hierarch}

We employ mixed-integer linear programming (MILP) to address the complexities of mixed bus fleet management, capturing various factors and constraints and providing a comprehensive and flexible framework for optimizing the integration of multiple kinds of buses in a fleet. Sivagnanam et al. \cite{sivagnanam2021minimizing} try to address the version of the problem without dynamic charging by minimizing the cost for all of the transit trips and non-service trips. Even state-of-the-art approaches do not scale for more than a few buses and a few transit trips (or transit blocks) per day when solving this MILP. To address this issue, we divide the problem into parts, by trying to assign buses to blocks hierarchically and by keeping the objective fairly simple. Our hierarchical approach enhances computational efficiency and scalability. In practice, for an NP-hard problem like ours, it's more feasible to tackle smaller sub-problems of size $n$ rather than a larger, more complex problem of size $m$, especially when $n << m$.  We generally observe that transit vehicles are of multiple types, and each of those types has a different operating efficiency. The hierarchical solution hinges on the fact that increasing the use of more efficient buses (like electric buses) decreases the overall cost of operations, as they cost less per mile. The hierarchy is ordered using $v^{op}$, and we use the operating efficiency threshold $v^{op}_{th}$ to divide the problem as follows:

\begin{enumerate}[leftmargin=*]
    \item \label{l1} \textbf{Level 1} We first attempt to assign the maximum blocks using buses which have $v^{op} \geq v^{op}_{th}$ (usually, electric buses).
    \item \label{l2} \textbf{Level 2} Buses that have $v^{op} < v^{op}_{th}$ (usually hybrid and diesel buses) are used next to serve the remaining blocks. 
\end{enumerate}

We could add more levels to the hierarchy based on multiple operating efficiency thresholds (for example, when using many types of electric buses with high variance in operating efficiency), but we refrain from doing so in our scenario. 

\noindent \textbf{Variables} We use two sets of binary variables and five sets of continuous variables. The first binary set, $a_{v, t}=1$ indicates that block $t$ is assigned to bus $v$, and $a_{v,(cp, s)}=1$ indicates electric bus $v$ is assigned to charging pole $cp$ at slot $s$. The second binary set, $m_{v, x_{1}, x_{2}}=1$ indicates that bus $v$ takes the non-service trip between a pair of transit blocks or charging slots, $x_{1}$ and $x_{2}$. We have five sets of continuous variables, $c_{s}^{v} \in\left[0, {\mathcal{M}_{v}}^{max}\right]$ represents the amount of energy charged to electric bus $v$ in time slot $s$. The second continuous variable set denotes the state of charge (SoC) of buses, $e_{s}^{v} \in\left[{\mathcal{M}_{v}}^{min}, {\mathcal{M}_{v}}^{max}\right]$ which is the battery level for electric buses and the fuel remaining for diesel buses. The continuous variable set, $g_s$ represents the recharging cost for the slot $s$. The two continuous variables $\delta_v$ and $u_t$ represent the distance covered by each bus, and the slack for each block assignment (the number of blocks unserved) respectively.

\subsection{Level 1}

We solve Level 1 of the problem using the following constraints.


\noindent \textbf{Constraints} Every transit block is served by exactly one bus:
\small
\begin{equation} \label{eq:assign}
\forall t \in \mathcal{T}: \sum_{v \in \mathcal{V}} a_{v, t} + u_t = 1
\end{equation}
\normalsize

Each bus should be used to its full potential when trying to serve the blocks, where $t^d$ is the transit block distance: 
\small
\begin{equation} \label{eq:delta}
\forall v \in \mathcal{V}: \delta_v = \sum_{t \in \mathcal{T}} a_{v, t} \cdot t^d
\end{equation}
\normalsize

Each charging slot is assigned at most one electric vehicle:
\small
\begin{equation}\label{eq:charge_v}
\forall(c p, s) \in \mathcal{C} : \sum_{\forall v \in \mathcal{V}: M_{v} \in \mathcal{M}^{\text {elec }}} a_{v,(cp, s)}  \leq 1
\end{equation}
\normalsize

 To ensure the converse of Eq~\ref{eq:charge_v}, that, one vehicle is not assigned to the multiple charging poles in the same slot:
\small
\begin{multline}\label{eq:charge_cp}
\forall v \in \mathcal{V}: \mathcal{M}_{v}\in\mathcal{M}^{\text {elec }}, \forall s \in S: \\
\sum_{\forall cp \in \mathcal{CP}} a_{v,(c p, s)} \leq 1
\end{multline}
\normalsize




We use $F\left(x_{1}, x_{2}\right)$ as described above in Section~\ref{sec:feasibility} to check for feasibility involving the movement between two blocks $x_1, x_2$ . The following constraint represents block feasibility:
\small
\begin{multline}\label{eq:feas}
\forall v \in \mathcal{V}, \forall x_{1}, x_{2} \in \mathcal{T} \cup \mathcal{C}, \neg F\left(x_{1}, x_{2}\right): a_{v, x_{1}}+a_{v, x_{2}} \leq 1
\end{multline}
\normalsize

When a bus $v$ is assigned to a pair of blocks or charging slots $x_1$ and $x_2$, and if $x_1$ and $x_2$ are consecutive assignments, then bus $v$ needs to take a non-service trip:
\small
\begin{multline}\label{eq:m1}
\forall v \in \mathcal{V}, \forall x_1, x_2 \in \mathcal{T} \cup \mathcal{C}, F(x_1, x_2): \\
m_{v, x_{1}, x_{2}} \geq a_{v, x_{1}}+a_{v, x_{2}}-1-\sum_{x \in \mathcal{T} \cup \mathcal{C}: x_{1}^{\text {start }} \leq x^{\text {start}}, x^{\text {end}} \leq x_{2}^{\text {start }}} a_{v, x}
\end{multline}
\begin{equation}\label{eq:m2}
    m_{v,x_1,x_2} \leq a_{v,x_1}, \quad m_{v,x_1,x_2} \leq a_{v,x_2}
\end{equation}
\normalsize

Equation \eqref{eq:m2} is used to limit the non-service trips to occur only between the blocks or charging slots the buses are assigned to.

The amount of energy charged to the vehicle $v$, at charging pole $cp$, is at a uniform rate of charge $c^v_s$ throughout the period of charging. We also ensure that the battery levels of electric buses remain between $\mathcal{M}_{v}^{min}$ and $\mathcal{M}_{v}^{max}$. Thus, for each electric bus $v$ in a slot $s$, the energy charged, $c_{s}^{v}$ is,
\small
\begin{multline}\label{eq:charging}
\forall v \in \mathcal{V}, s \in \mathcal{S}:
    c_{s}^{v} \leq \sum_{(cp, s) \in \mathcal{C}} q^{max}_{cp} \cdot a_{v,(cp, s)}
\end{multline}
\normalsize

Next, we include the cost of recharging the vehicles, which can be affected by the dynamic electricity prices:
\small
\begin{equation}\label{eq:g_s}
    \forall s \in S, g_s = w_s \cdot \sum_{v\in \mathcal{V}} c^v_s
\end{equation}
\normalsize
For all diesel vehicles, 
\small
$\mathcal{M}_v \in \mathcal{M}^{diesel}$, $c^v_{s_n} = 0$.

Since $e_{s}^{v} \in\left[{\mathcal{M}_{v}}^{min}, {\mathcal{M}_{v}}^{max}\right]$, it is ensured the bus never runs out of charge and recall, that we represent the amount of energy used to move between $x_1, x_2$ by $E\left(v, \mathcal{T}\left(x_{1}, x_{2}\right)\right)$.
For the first time slot, $n = 0$, $e_{0}^{v} =  {\mathcal{M}_v}^{max}$. For any other $n^{th}$ time slot $s_{n}$, and for a bus $v$, we can find the amount of energy remaining $e_{s_{n}}^{v}$, as
\small
\begin{multline}\label{eq:energy}
\forall v \in \mathcal{V}, s \in \mathcal{S}:\\
e_{s_{n}}^{v}= e_{s_{n-1}}^{v} + c_{s_{n}}^{v} - \sum_{t \in \mathcal{T}: s_{n}^{\text {start }}<t^{\text {end }} \leq s_{n}^{\text {end }}} a_{v, t} \cdot E(v, t) \\
-\sum_{x_{1}, x_{2}: s_{n}^{\text {start }}<x_{2}^{\text {start }} \leq s_{n}^{\text {end }}} m_{v, x_{1}, x_{2}} \cdot E\left(v, \mathcal{T}\left(x_{1}, x_{2}\right)\right)
\end{multline}
\normalsize

\normalsize

To find the cost of replenishing the buses, we find the energy left at the end of all transit blocks (at the last slot of the day, $s_{last} $) and replenish each vehicle to its maximum capacity:
\small
\begin{multline}\label{eq:gs_end}
g_{s_{last}} = \sum_{v\in \mathcal{V}:\mathcal{M}_v \in \mathcal{M}^{elec}}({\mathcal{M}_{v}}^{max}-e^v_{s_{last}}) \cdot w_{s_{last}} \\
+ \sum_{v\in \mathcal{V}:\mathcal{M}_v \in \mathcal{M}^{diesel}} ({\mathcal{M}_{v}}^{max} - e^v_{s_{last}}) \cdot K^{diesel}
\end{multline}
\normalsize
where $K^{diesel}$ is the cost of diesel per gallon.

To ensure seating capacity constraints on a block are met:
\small
\begin{equation}\label{eq:seating}
    v^s \geq t^s
\end{equation}
\normalsize

\noindent \textbf{Objective} We maximize the distance covered by the more efficient buses, maximizing $\delta_{v}$ the distance covered by them.
\small
\begin{equation}\label{eq:obj1}
\max \sum_{v \in \mathcal{V}} \delta_v
\end{equation}
\normalsize

\subsection{Level 2}

In this level, we assign buses below $v^{op}_{th}$ to all the transit blocks that remain after Level 1 assignments are done, while minimizing energy use. Here, we can reduce the variables and constraints needed, and the problem becomes simpler.

\noindent \textbf{Variables} We use all variables except two, the slack variable $u_t$  as we assign all blocks to a bus, and $\delta_v$ which was used to measure distance traveled by a bus.

\noindent \textbf{Constraints} We modify Eq.~\eqref{eq:assign} to assign buses to all blocks and ensure feasibility:
\small
\begin{equation} \label{eq:assign_mod}
\forall t \in \mathcal{T}: \sum_{v \in \mathcal{V}} a_{v, t} = 1
\end{equation}
\normalsize

Eq~\ref{eq:delta} is no longer required as we do not need to maximize the distance. We use the other constraints \eqref{eq:charge_v} to \eqref{eq:seating}.

\noindent \textbf{Objective}
We modify the objective completely to align with reducing the operating costs. Thus, we minimize the total replenishing cost.
\small
\begin{equation}\label{eq:obj2}
\min \sum_{s \in \mathcal{S}} g_s
\end{equation}
\normalsize
\subsection{Maintaining Feasibility}

We also need to check and maintain the feasibility of the solution using the hierarchical approach. A solution set could become infeasible when solved using strict hierarchical rules. For example, say we have a problem set with 4 transit blocks ($t_1, t_2, t_3)$, all starting at the same time, and $t_4$ starting just after the end of the previous transit blocks and from the same location. The fleet is made up of 2 electric buses and 1 diesel bus. Say, $t_1, t_4$ are $\Delta_1$ in length, while the other transit blocks are $\Delta_2$ in length, where $\Delta_1 > \Delta_2$. Also, the electric buses cannot serve $t_1$ and $t_4$ consecutively due to battery constraints. Now, in the hierarchical solution, since we maximize the distance covered by electric buses, they will be assigned to $t_1$ and $t_4$, leaving $t_2, t_3$. 1 diesel bus cannot serve both blocks at the same time, making it infeasible. Whereas, if all the buses and blocks are considered together and solved non-hierarchically, we find a feasible solution where electric buses serve $t_1, t_3$ and $t_2, t_4$ are served by the diesel bus. Again, these issues can occur only in very specific scenarios and are not common, as we show in our experiments.

\begin{algorithm}
\caption{Hierarchical Iteration Algorithm}\label{algo:iteration_h}
\begin{algorithmic}[1]
\Function{Hierarchical Iteration}{$\mathcal{T}, \mathcal{V}, \mathcal{C}$}
    \State $feas \gets True$ 
    \State $\mathcal{V}_1 \gets \{v \in \mathcal{V}:\mathcal{M}_v \in \mathcal{M}^{elec}\}$
    \State $\mathcal{V}_2 \gets \{v \in \mathcal{V}:\mathcal{M}_v \in \mathcal{M}^{diesel}\}$
    \While {$feas$}
        \State $\mathcal{A}_1 \gets$ Level 1({$\mathcal{T}, {\mathcal{V}_1}, \mathcal{C}$}); 
        $\mathcal{T}' \gets \mathcal{T} \setminus \{t \in \mathcal{A}_1\}$
        \State $\mathcal{A}_2 \gets$ Level 2({$\mathcal{T}', \mathcal{V}_2, \mathcal{C}$})
        
        \State $\mathcal{V}_1 \gets \mathcal{V}_1 \setminus \{\hat{v}\}$; $\mathcal{V}_2 \gets \mathcal{V}_2 \cup \{\hat{v}\}$
        
        \State $feas \gets feasible(\mathcal{A}_2)$
        
    \EndWhile
    \State \Return $\mathcal{A}_1 \cup \mathcal{A}_2$
    \EndFunction
\end{algorithmic}
\end{algorithm}

To achieve feasibility, we introduce an iterative process to the solution method. We solve the problem using the hierarchical method described above and then check the feasibility of the solution. If feasible it is accepted. If not feasible, we remove one of the least efficient buses from Level 1 and add it to Level 2. This process continues till we get a feasible solution. If all buses are removed from Level 1 to Level 2, the problem becomes a non-hierarchical MILP, where feasibility is guaranteed, provided there are enough buses to serve the blocks. Algorithm~\ref{algo:iteration_h} shows this process.

\subsection{Warm start using initial solution}\label{sec:greedy}

To speed up the Hierarchical solution we form a Greedy approach that can be used to provide an initial set of solutions for warm starting. The greedy solution can provide a good approximation to the optimal solution, giving the solver a good solution to start with which can help in finding or get near the optimal solution quicker.
\begin{algorithm}
\caption{CanServe Algorithm}
\label{algo:can_serve}
\begin{algorithmic}[1]
\Function{CanServe}{$\mathcal{A}, \mathcal{T}, \mathcal{C}, v, x$}
\If{$\{ {x} \in \mathcal{T \cup C} | \langle v, {x} \rangle \in \mathcal{A} \}$}
    \State \textit{Previous} $\gets \{ \hat{x} \in \mathcal{T \cup C} | \langle v, \hat{x} \rangle \in \mathcal{A} \wedge \hat{x}^{end} \leq x^{start} \}$ 
    
    \If{\textit{Previous} $\neq \emptyset$}
        \State $x_{prev} = argmax_{\hat{x} \in Previous} \hat{x}^{end}$
        \If{$x_{prev}^{end} + D(x_1, x_2) \ge x^{start}$}
            \State \Return False
        \Else
            \State \Return $e^{v}_s - E(v, t) \ge {\mathcal{M}_v}^{min}$
        \EndIf
    \EndIf
\EndIf
\EndFunction
\end{algorithmic}
\end{algorithm}

\textbf{CanServe Algorithm}: The Algorithm~\ref{algo:can_serve} determines if a given bus can serve a particular block. It checks if the bus can serve the current transit block $t$ based on its previous assignments while ensuring it has enough energy to serve the current block $t$. The computational complexity of the CanServe Algorithm is O(1).

\begin{algorithm} 
\caption{BiasedCost Algorithm}\label{algo:biased}
\begin{algorithmic}[1]
\Function{BiasedCost}{$\mathcal{A}, \mathcal{T}, \mathcal{C}, v, x, wp$}
\State \textit{Previous} $\gets \{ \hat{x} \in \mathcal{T \cup C} | \langle v, \hat{x} \rangle \in \mathcal{A} \wedge \hat{x}^{end} \leq x^{start} \}$ 
\State \textit{After} $\gets \{ \hat{x} \in \mathcal{T \cup C} | \langle v, \hat{x} \rangle \in \mathcal{A} \wedge \hat{x}^{start} \geq x^{end} \}$ 

\If{\textit{Previous} $\neq \emptyset$}
\State $x_{prev} = argmax_{\hat{x} \in Previous} \hat{x}^{end}$
\State $m_{prev} = \mathcal{T}( x^{destination}_{prev}, x^{origin})$
\EndIf
\If{\textit{After} $\neq \emptyset$}
\State $x_{prev} = argmax_{\hat{x} \in After} \hat{x}^{end}$
\State $m_{prev} = \mathcal{T}(  x^{destination}, x^{origin}_{prev})$
\EndIf
\State $g \gets E(v, t) + E(v, m_{prev}) + wp \cdot (x^{start} - x^{end}_{prev})$
\State \Return $g$
\EndFunction
\end{algorithmic}
\end{algorithm}

\textbf{BiasedCost Algorithm}: Algorithm~\ref{algo:biased} calculates the biased cost of assigning a block to a bus based on various factors. The algorithm considers the energy required to serve the transit block and the cost of a non-service trip from a previous block. The movement time is modified by the Wait Penalty, $wp$, which helps to prioritize assigning buses that have to travel a shorter time to get to and from a block. The biased cost is the sum of the energy used for transit block, non-service trip, and the wait penalty ($wp$). The complexity of the BiasedCost Algorithm is O(1).

\begin{algorithm}
\caption{Greedy Assignment Algorithm}\label{algo:greedy}
\begin{algorithmic}[1]
\Function{GreedyAssignment}{$\mathcal{T}, \mathcal{V}, \mathcal{C}, wp$}
    \State $\mathcal{A} \gets \emptyset$
    \For{$t \in \mathcal{T}$}
        \State $\mathcal{F} \gets \{v \in \mathcal{V} | \textbf{CanServe}(\mathcal{A}, \mathcal{T}, \mathcal{C}, v, t)\}$
        \State $\hat{v} \gets \arg\min_{v \in \mathcal{F}} \textbf{BiasedCost}({\mathcal{A}, \mathcal{T}, \mathcal{C}, v, t, wp})$
        
        \If{$\hat{v} \neq \emptyset$}
            \State $\mathcal{A} \gets \hat{v}$
            \If{$ \hat{v} \in \mathcal{M}^{elec}\ \&\ e^{\hat{v}}_s - E(\hat{v}, t) \le \mathcal{V}^{th}_{\mathcal{M}_v}$}
                \For{$(cp, s) \in C$}
                    \If{$s \ge t^{end}$}
                        \State $e^{v}_s = e^{v}_s + q^{max}_{cp} \cdot (s^{end} - s^{start})$
                    \EndIf
                \EndFor
            \EndIf
        \EndIf
    \EndFor

    \EndFunction
\end{algorithmic}
\end{algorithm}

\textbf{Greedy Algorithm}: The Greedy Assignment Algorithm~\ref{algo:greedy} efficiently assigns transit blocks to the available buses in the fleet, considering the energy consumption and charging requirements. We provide it with the fleet composition ($\mathcal{V}$), the set of charging slots ($\mathcal{C}$), and charging rates. The algorithm iteratively processes each block in the set $\mathcal{T}$, and for each block, it first determines the set of feasible buses $\mathcal{F}$, found by using the \textbf{CanServe} Algorithm. \textbf{BiasedCost} is utilized to identify the bus with the minimum biased cost, $\hat{v}$, and is used for serving the block.
If the selected bus is an electric bus and its SoC remaining after serving the block is less than the minimum SoC,  $\mathcal{V}^{th}_{\mathcal{M}_v}$ a charging slot is assigned from the set $C$.
The computational complexity of the Greedy Assignment Algorithm is O($|\mathcal{T}|\cdot |\mathcal{V}| \cdot |\mathcal{C}|$), where $|\mathcal{T}|$ is the set of transit blocks,  $|\mathcal{V}|$ is the bus fleet and |$\mathcal{C}$| is the set of charging slots. The algorithm is efficient in solving large-scale problems, very quickly.

\section{Experiments}

\subsection{Data Preparation}

We use the transit data from our partner agency, Chattanooga Area Regional Transportation Authority (CARTA) to perform the experiments. We use their General Transit Feed Specification (GTFS) to obtain the transit blocks and their fleet composition. GTFS is a widely adopted data format used by transit agencies to define and share public transportation schedules and related geographic information. It provides comprehensive information about transit routes, transit blocks, stops, trips, and schedules. 

For testing, we use one month of their recent deployment data from February 2024. During this time, CARTA managed an average of 430 trips, grouped into 42 blocks, utilizing a mixed fleet of 31 diesel buses, and 4 electric buses (BYD K9M model). It is important to note that not all buses were available every day due to maintenance downtime. The diesel buses averaged a fuel consumption of $2.53$ miles/gal, whereas the electric buses averaged $0.56$ miles/kWh. We use the U.S. Energy Information Administration (EIA) \footnote{\hyperlink{U.S. EIA}{https://www.eia.gov/energyexplained/units-and-calculators/energy-conversion-calculators.php}} estimated conversion rate of 1 gal of diesel $=$ 37.1 kWh to compare and compute the usage of electric and diesel bus. We set diesel costs at $\$4.2$/gal and utilized CARTA's actual bus usage data for each day of February for our simulations. The BYD K9M electric buses in CARTA's fleet have a battery capacity of 310 kWh, and the 2 chargers are located at the central depot, with a maximum charging rate of 80 kWh.  

Since Chattanooga uses a uniform cost of electricity, we set the electricity price at $\$0.12795$/kWh. For other experiments, we use a time-of-use price, where the peak hours (6am to 10pm) are at $\$0.14660$/kWh, and the off-peak hours (time outside peak hours) are priced at $\$0.12795$/kWh. We set $v^{op}_{th} = \min{(v^{op})}\ \forall v\in \mathcal{V}:\mathcal{M}_v \in \mathcal{M}^{elec}$, i.e. the minimum operating efficiency of the electric bus fleet. The minimum seat requirement for every block is set at 30, and all buses in CARTA's fleet have a minimum of 30 seats.

 \begin{figure}[htp]
    \centering
    \includegraphics[width=0.49\textwidth]{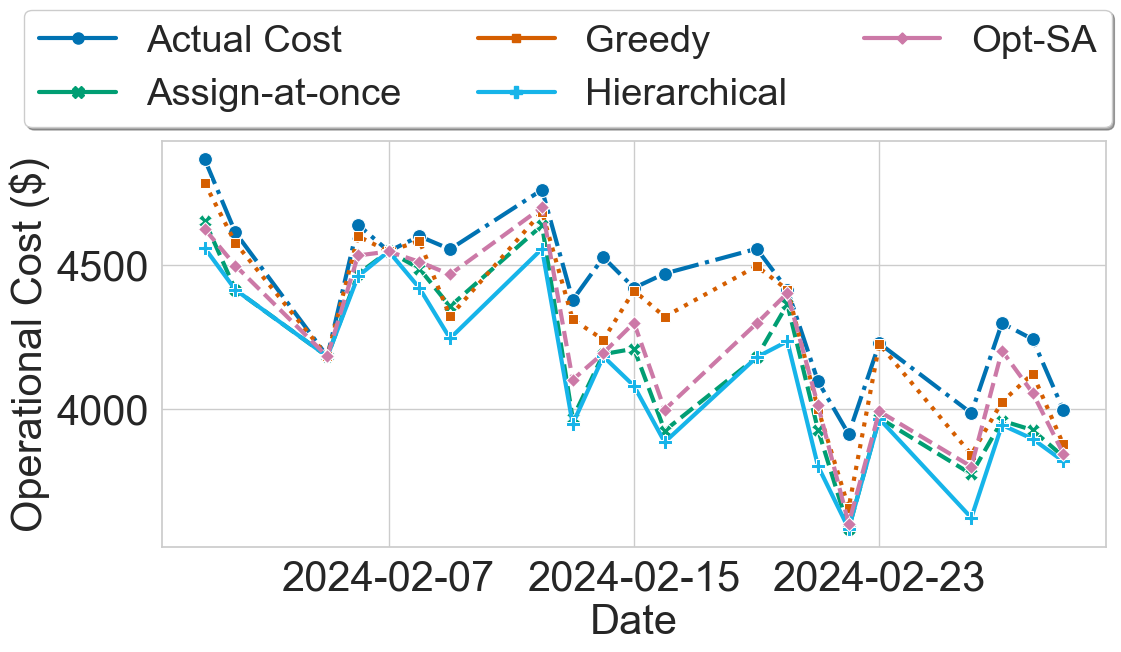}
    \caption{Comparison of all methods across Feb 2024 for actual block assignments used by CARTA. With a maximum fleet size of 4 electric and 31 diesel buses, which can vary daily. The electric buses have a battery capacity of 310 kWh, and 2 chargers with a maximum charging rate of 80 kWh. We compare the operational cost for the Hierarchical and the baseline models. Hierarchical provides the best operational cost for all days}
    \label{fig:all_comp}
\end{figure}

 \begin{figure}[htp]
    \centering
    \includegraphics[width=0.44\textwidth, height=3cm]{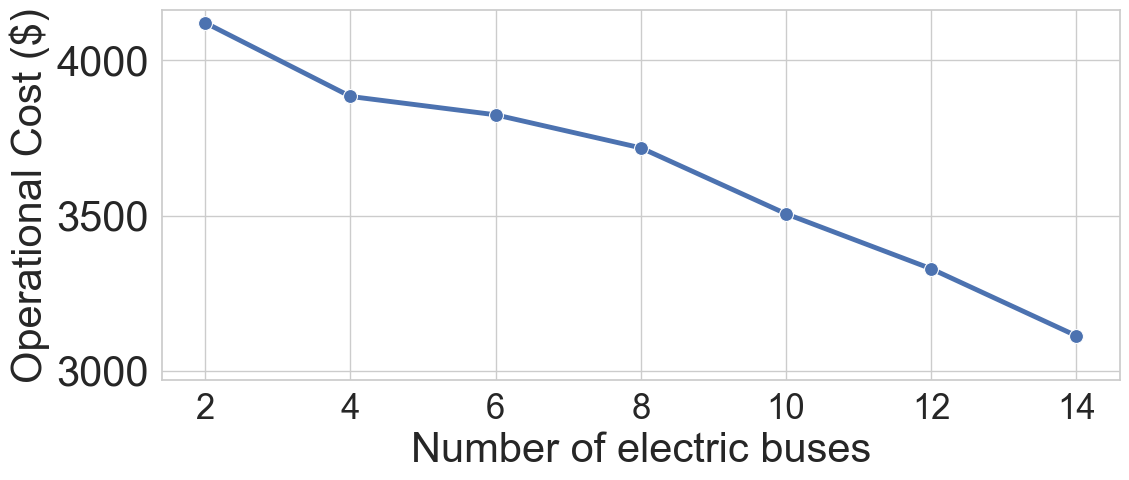}
    \caption{Effect of changing the number of electric buses for 2024-02-16, which has 47 blocks. The number of diesel buses are decreased proportionately to maintain a total available fleet size of 27}
    \label{fig:elec_buses}
\end{figure}

\subsection{Baselines}

We refer to our proposed algorithm as \textbf{Hierarchical}, and use the following baselines for comparison.

\noindent\textbf{Actual cost}: The actual operational cost incurred by CARTA for the day's service, based on the buses deployed to specific transit blocks, using the average efficiency metric of the buses as described above. 

\noindent \textbf{Greedy}: The Greedy assignment as described in section \ref{sec:greedy} assigns buses based on  a first-come-first-serve basis, prioritizing electric buses over diesel buses. we use waiting penalty, $wp = 0.001$, and the SoC threshold $\mathcal{V}^{th}_{\mathcal{M}_v} = 0.5 \cdot \mathcal{M}^{min}_v \ \forall v \in \mathcal{V}: M_{v} \in \mathcal{M}^{\text {elec }}$. We choose half of the maximum capacity to be a threshold for charging as a precaution, as we do not want the buses to go below their minimum allowable charge ${\mathcal{M}}^{min}_v$.

\noindent \textbf{Assign-at-once}: Solving only using Level 2 of the hierarchical solution described in 
Section~\ref{sec:hierarch}. It is a one-step solution to the problem, proposed above. We use the Greedy solution as an initial solution to warm start it - which helps improve the solution search time.

\noindent \textbf{Opt-SA}: This state-of-the-art baseline is the simulated annealing solution implemented by Sivagnanam et. al \cite{sivagnanam2021minimizing}, which uses a greedy initial solution and follows the simulated annealing method to swap blocks to reduce operational costs. As used by the author, an iteration limit of 40000 for the algorithm, along with the other parameters.

\subsection{Multiple Days Analysis}

We use an off-the-shelf MILP solver --- CPLEX \cite{cplex2009v12} to solve all the MILPs. The MILP solving times are limited to a maximum time of 15 minutes. 15 minutes proved to be sufficient for use, as we obtained an average relative gap of $0.1\%$. Using Fig.~\ref{fig:all_comp} we analyze bus deployment from Feb 1 to Feb 29, 2024, and evaluate Actual Cost, Greedy, Opt-SA, Assign-at-once, and Hierarchical solution methods. The Hierarchical method produces the least operational cost while maintaining feasible solutions. On average, the proposed hierarchical solution is better than the actual cost of CARTA by $\mathbf{6.25\%}$, and $\mathbf{2.65\%}$ better than Opt-SA. In the best case, it outperforms Opt-SA by $\mathbf{8.37\%}$, and the actual cost of CARTA by $\mathbf{15.37\%}$.


\subsection{Impact of Increasing Electric Vehicles (EVs)}
With the increasing electrification of vehicles, we attempt to estimate the possible savings we can get by having electric buses form a larger portion of the transit fleet. We measure the impact of increasing electric bus numbers on operational costs in Fig.~\ref{fig:elec_buses} on a day's assignment (Feb 16, 2024) as it has operational costs close to the mean operational cost for the month. As we can see, the operational costs are reduced significantly when using electric buses and significantly reduce tailpipe emissions. Compared to a fleet with 2 electric buses, a fleet with 14 electric buses saves $\$1009.08$ in operational cost, and $25.64$ metric tons of CO$_2$ emissions daily according to U.S. Environmental Protection Agency (EPA) Greenhouse Gas Equivalencies.\footnote{\hyperlink{U.S. EPA}{https://www.epa.gov/energy/greenhouse-gases-equivalencies-calculator-calculations-and-references}}

 \begin{figure}[htp]
    \centering
    \includegraphics[width=0.4\textwidth, height=3.1cm]{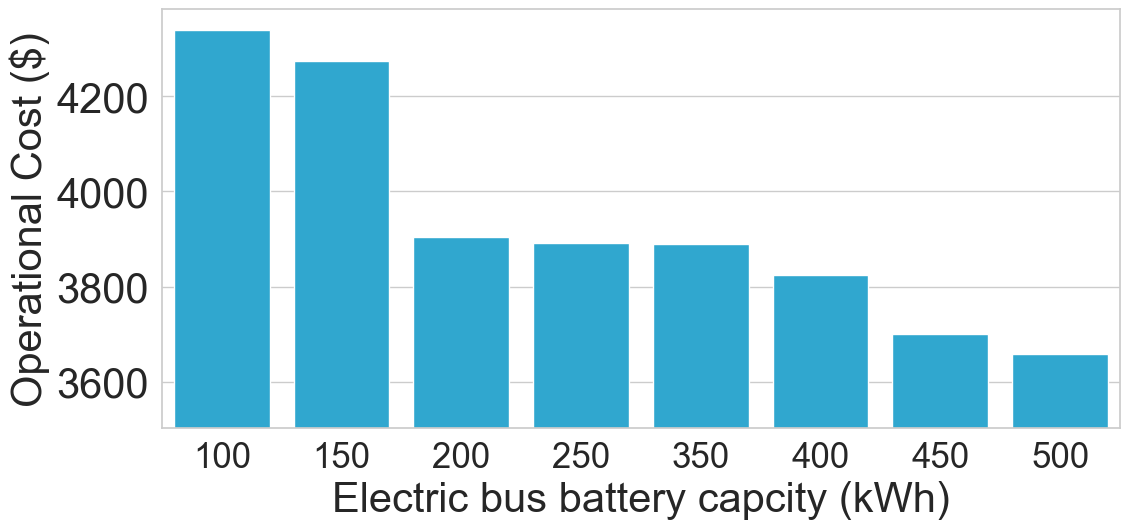}
    \caption{Effect of changing the electric bus battery capacity on the operational cost}
    \label{fig:batt_cap_change}
\end{figure}

\vspace{-5mm}
 \begin{figure}[htp]
    \centering
    \includegraphics[width=0.4\textwidth, height=3.1cm]{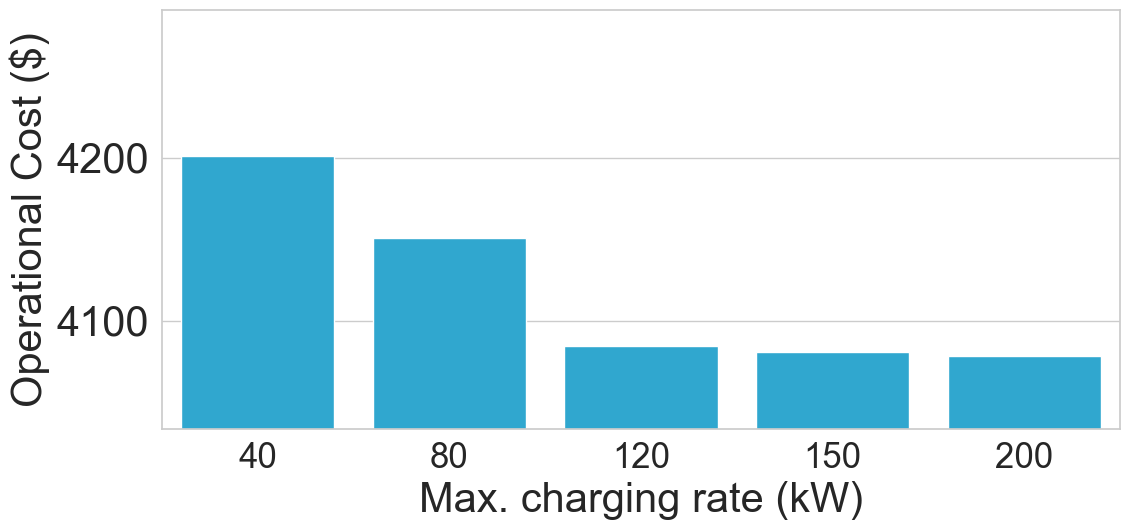}
    \caption{Effect of changing the maximum charging rate of the chargers on the operational cost}
    \label{fig:maxq}
\end{figure}

\vspace{-5mm}
 \begin{figure}[htp]
    \centering
    \includegraphics[width=0.44\textwidth, height=3.1cm]{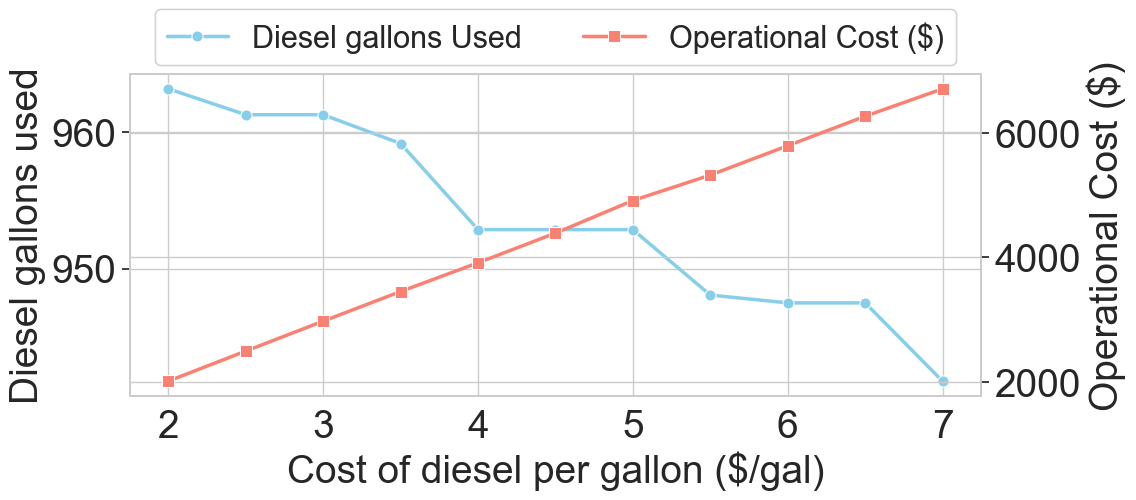}
    \caption{Effect of changing the diesel cost per gallon on the diesel bus usage (in terms of gallons used). The operational cost increases due to the increase in diesel cost.}
    \label{fig:diesel_cost}
\end{figure}

\subsection{Sensitivity Analysis}
\begin{figure}[htp]
    \centering
    \includegraphics[width=0.44\textwidth, height = 3.4cm]{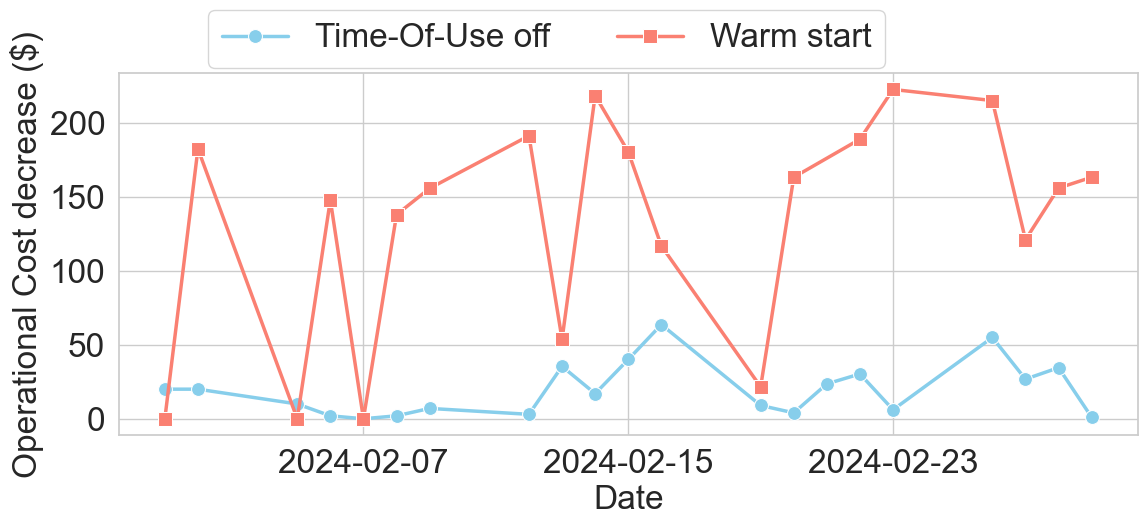}
    \caption{Operational cost decreases when using Greedy solution for warm start against when it is not used. Similarly, we show the difference when the Time-of-Use price is used compared to when it is not used.}
    \label{fig:tou_ws}
\end{figure}

\begin{figure}[htp]
    \centering
    \includegraphics[width=0.38\textwidth, height = 2.5cm]{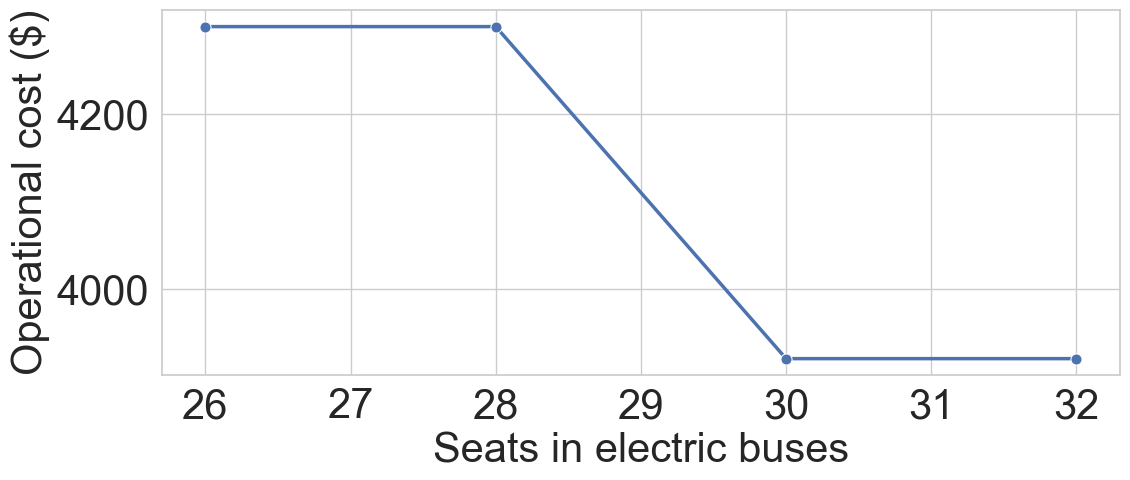}
    \caption{Effect of changing the seats in electric buses for 2024-02-16, which has a minimum requirement of 30 seats on all blocks}
    \label{fig:seats}
\end{figure}

We test the sensitivity and adaptability of the Hierarchical model to changes in key parameters such as battery limits in Fig.~\ref{fig:batt_cap_change}, maximum charging rate in Fig.~\ref{fig:maxq}, Fig.~\ref{fig:seats} for seats available and Fig.~\ref{fig:diesel_cost} for diesel price changes. All the experiments for sensitivity analysis are done on CARTA's actual assignments for Feb 16, 2024,  which has 47 blocks and uses 4 electric and 23 diesel buses and has an operational cost close to the mean operational cost for the month. In Fig.~\ref{fig:tou_ws} we also observe that performing a warm start (using the Greedy solution) reduces the operational cost, given the same time limits. The use of Time-of-Use price causes a minor increase in the operational cost, as the electric buses will try to recharge during the peak hours, to cover more distance, incurring a higher electricity cost.

\vspace{-5mm}
\section{Conclusion}\label{sec:conc}



We present a novel formulation for optimizing public transit fleets on fixed lines under dynamic grid pricing and seating constraints, using a hierarchical approach. This method demonstrates superior performance over single-stage solutions and heuristic-based methods like simulated annealing by strategically prioritizing higher operational and fuel efficiency buses and adopting a two-step (or multi-step) resolution process. We also leverage the use of a Greedy initial solution to enhance solution quality within limited solving time. Our findings demonstrate how transit agencies can alleviate the complex challenge of managing mixed fleet transit infrastructure, and better manage their operations. We provide a daily average saving in operational costs of $\$275.24$ and an estimated reduction of 3.05 metric tons of CO$_2$ in tailpipe emissions every day (according to U.S. EPA greenhouse equivalents) compared to existing schedules used by our partner transit agency CARTA.


\section{Acknowledgements}

This material is based upon work sponsored by the National Science Foundation (NSF) under grant CNS-1952011, the Department of Energy under Award Number DE-EE0009212, and Chattanooga Area Regional Transportation Authority (CARTA). Results presented in this paper were obtained using the Chameleon Testbed supported by the NSF.

\bibliographystyle{IEEEtran}
\bibliography{main}

\end{document}